\documentclass[12pt,a4paper]{amsart}
\usepackage[latin1]{inputenc}
\usepackage{latexsym}
\usepackage{color,graphicx,shortvrb}
\usepackage{amsmath, amssymb}
\usepackage{amsfonts}
\usepackage[colorlinks, bookmarks=true]{hyperref}

\newtheorem{theorem}{Theorem}[section]
\newtheorem{lemma}[theorem]{Lemma}
\newtheorem{proposition}[theorem]{Proposition}
\newtheorem{corollary}[theorem]{Corollary}

\theoremstyle{definition}
\newtheorem{definition}[theorem]{Definition}

\newcommand{\rank}{\mathbf{rank}}

\title{On  strict inclusion relations between approximation and interpolation spaces}
\author{J. M. Almira}

\begin{document}

 \baselineskip=16pt

\numberwithin{equation}{section}

\maketitle \markboth{Strict inclusion relations}{J. M. Almira}

\begin{abstract}
We study strict inclusion relations between approximation and interpolation spaces.
 \end{abstract}

\section{Introduction}
Suppose $X$ is a quasi-Banach space, and let
$A_0\subset A_1 \subset \ldots \subset X$
be an infinite chain of subsets of $X$, where all inclusions are
strict. We say that $(X,\{A_n\})$ is an {\it approximation scheme}
(or that $(A_n)$ is an approximation scheme in $X$) if:
\begin{itemize}
\item[$(i)$] There exists a map $K:\mathbb{N}\to\mathbb{N}$ such that $K(n)\geq n$ and $A_n+A_n\subseteq A_{K(n)}$ for all $n\in\mathbb{N}$
(we can assume that $K$ is increasing).

\item[$(ii)$] $\lambda A_n\subset A_n$ for all $n\in\mathbb{N}$ and all scalars $\lambda$.

\item[$(iii)$] $\bigcup_{n\in\mathbb{N}}A_n$ is dense in $X$.
\end{itemize}
Approximation schemes were introduced  by Butzer and Scherer \cite{butzer_scherer} in 1968 and, independently,  by Y. Brudnyi and N. Kruglyak under the name of ``approximation families'' in 1978 \cite{brukru}, and
popularized by Pietsch in his seminal paper \cite{pietsch}, where the approximation spaces  $A_p^r(X,{A_n})=\{x\in X: \|x\|_{A_p^r}=\|\{E(x,A_n)\}_{n=0}^\infty\|_{\ell_{p,r}}<\infty\}$ were studied. Here, $$\ell_{p,r}=\{\{a_n\}\in\ell_\infty: \|\{a_n\}\|_{p,r}=\left[\sum_{n=1}^\infty n^{rp-1}(a_n^*)^p\right]^\frac{1}{p}<\infty\}$$ denotes the so called Lorentz sequence space. Pietsch was interested in the parallelism that exists between the theories of approximation spaces and interpolation spaces (e.g., he proved embedding, reiteration and representation results for approximation spaces). In all these cases the authors imposed condition $A_n+A_m\subseteq A_{n+m}$, which implies $K(n)=2n$. Simultaneously and also independently, Ti\c{t}a \cite{tita0} studied, from 1971 on, for the case of approximation of linear operators by finite rank operators,  a similar concept, based on the use of symmetric norming functions $\Phi$ 
and the sequence spaces defined by them, $S_{\Phi}=\{\{a_n\}:\exists \lim_{n\to\infty}\Phi(a_1^*,a_2^*,\cdots,a_n^*,0,0,\cdots)\}$.
The concept of approximation scheme given in the present paper was introduced by Almira and Luther \cite{almiraluther2}, \cite{almiraluther1} a few years ago. They also created a  theory for generalized approximation spaces via the use of general sequence spaces $S$ (that they named ``admissible sequence spaces'') and the definition of the approximation spaces $A(X,S,\{A_n\})=\{x\in X: \|x\|_{A(X,S)}=\|\{E(x,A_n)\}\|_S<\infty\}$. Other papers with a similar spirit of generality have been written by Ti\c{t}a \cite{tita_cluj_99} and Pustylnik \cite{pustylnik}, \cite{pustylnik2}.  It is important to remark that, due to the importance of the so called direct and inverse theorems in approximation theory (also named ``central theorems in approximation theory''), the idea of defining approximation spaces is a quite natural one. This has produced the negative effect that many unrelated people has though on the same things at different places and different times, and some papers in this subject partially overlap.

A fundamental part of the  theory developed by the authors of the above mentioned papers consists of the study of the embeddings between the involved spaces. Concretely, Pietsch proved that the embedding $A_p^r(X,{A_n})\hookrightarrow A_q^s(X,{A_n})$ holds true whenever $r>s>0$ or $r=s$ and $p<q$. This, in conjunction with the central theorems in approximation theory, which state a strong relation between smoothness of functions $f$ (compactness of operators $T$, respectively) and  fast decay of approximation errors $E(f,A_n)$ (approximation numbers $a_n(T)$, respectively), has been used to speak about the scale of smoothness (compactness, respectively) defined by an approximation scheme $(X,\{A_n\})$. Concretely, it is assumed that membership to the approximation space $A_p^r(X,\{A_n\})$ is a concept of smoothness (compactness if $X=B(Y_1,Y_2)$ and $A_n=\{T\in B(Y_1,Y_2): \rank (T)<n\}$). Thus, although in many concrete cases of approximation schemes $(X,\{A_n\})$ there exist some results, such as the
representation theorems or the characterizations in terms of moduli
of smoothness or in terms of interpolation spaces, that allow to prove
that all inclusions are strict (hence distinct parameters $p,r$ define
distinct smoothness concepts), the full truth says us that there is no
general result guaranteing this property. The main goal of this paper is to study  some results proving that, for certain classes of approximation schemes $(X,\{A_n\})$ and sequence spaces $S$, if $S_1\subset S_2\subset c_0$ (with strict inclusions) then $A(X,S_1,\{A_n\})$ is properly contained into $A(X,S_2,\{A_n\})$. It is clear that the same kind of questions (about strictness of the inclusions between distinct spaces), can be addressed for the context of interpolation spaces. We devote last section of this paper to initiate a study of this problem for Petree's real interpolation method.

Previous to any work in the above mentioned direction, it is necessary to solve the following question: Under which conditions on the approximation scheme $(X,\{A_n\})$ and the admissible sequence space $S$ we have that $A(X,S,\{A_n\})$ is properly contained into $X$? This question has been fully solved in the Banach setting by Almira and Oikhberg\cite{almira}. Concretely,  we can use Theorem 3.3 from \cite{almira} to guarantee that, if $X$ is Banach and $(X,\{A_n\})$ satisfies Shapiro's theorem, for all $(\varepsilon_n)$ decreasing to zero, there exists $x\in X$ such that $E(x,A_n)\geq \varepsilon_n$ for all $n=0,1,\cdots$. In particular, if $X$ is Banach and $S$ is an admissible sequence space properly contained into $c_0$, then $A(X,S,\{A_n\})$ is a proper subspace of $X$ whenever $(X,\{A_n\})$ satisfies Shapiro's theorem. Recall that $(X,\{A_n\})$ satisfies Shapiro's Theorem if and only if for all sequence $\{\varepsilon_n\}\in c_0$ there exists $x\in X$ such that $E(x,A_n)\not=\mathbf{O}(\varepsilon_n)$. These schemes were characterized in \cite[Corollary 3.6]{almira} as those verifying $E(S(X),A_n)=\sup_{\|x\|=1}E(x,A_n)=1$ for all $n$ (e.g., Riesz's Lemma implies that nontrivial linear approximation schemes satisfy Shapiro's theorem). Moreover, in the same paper the authors show many examples of schemes verifying Shapiro's theorem. In particular, the process of approximation by finite rank operators satisfies Shapiro's theorem \cite[Corollary 6.24]{almira}. Finally, we should also mention that \cite[Proposition 5.6]{almira} guarantees that, if $K(n)=cn$ for a certain finite constant $c>0$ and $(X,\{A_n\})$ satisfies Shapiro's theorem, then inclusion $A_p^r(X,{A_n})\hookrightarrow A_q^s(X,{A_n})$ is strict whenever $r>s>0$. We will present an easier proof of this fact for the special case of approximation schemes $(X,\{A_n\})$ verifying $E(S(X)\cap A_{n+1},A_n)>c>0$ for infinitely many natural numbers $n$ and a fixed constant $c>0$, a condition already used by Brudnyi in \cite{brundykrugljak}.

\section{Preliminary definitions and notations }

\begin{definition} Let $S$ be a real linear space of sequences $\{a_{n}\}_{n=0}^{\infty}%
\subseteq\mathbb{R}$ (with element-wise defined operations), equipped with a
quasi-norm $\Vert\,.\,\Vert_{S}\,.$ $S$ is called admissible sequence space
[for the approximation scheme $(X,\{A_{n}\})$] if the following
assumptions are satisfied:
\begin{itemize}
\item[$(A1)$] All finite sequences $\{a_{n}\}_{n=0}^{N}$ belong to $S$.

\item[$(A2)$]  If $0\leq a_{n}\leq b_{n}$ for all $n=0,1,2,\ldots$ and $%
\{b_{n}\}_{n=0}^{\infty}\in S$, then $\{a_{n}\}_{n=0}^{\infty}\in S$ and $%
\Vert\{a_{n}\}\Vert_{S}\leq\Vert\{b_{n}\}\Vert_{S}\,.$

\item[$(A3)$]  If $a_{0}\geq a_{1}\geq a_{2}\geq\ldots\geq0$ and $\{a_{K(n)}\}_{n=0}^{%
\infty}\in S$ then $\{a_{n}\}\in S$ and
\begin{equation*}
\Vert\{a_{n}\}_{n=0}^{\infty}\Vert_{S}\leq
C_{S}\,\Vert\{a_{K(n)}\}_{n=0}^{\infty}\Vert_{S}\,,
\end{equation*}
where $C_{S}$ is a constant, which only depends on $S$ and $%
\{K(n)\}_{n=0}^{\infty}\,.$
\end{itemize}
\end{definition}

\begin{definition} Let $S$ be an admissible sequence space for the
approximation scheme $(X,\{A_{n}\})$. Then the space
\begin{equation*}
A(X,S,\{A_n\}):=\{x\in X:\{E(x,A_{n})\}_{n=0}^{\infty}\in
S\}\,,\,\
\end{equation*}
endowed with$\,\Vert f\Vert_{A(X,S)}:=\Vert\{E(x,A_{n})\}_{n=0}^{%
\infty }\Vert_{S}\,$ is called (generalized) approximation space.
\end{definition}

A theory for these spaces was developed by Almira and Luther \cite{almiraluther2}. In particular, they proved the following properties:
 \begin{itemize}
 \item[$(P1)$] $(A(X,S,\{A_n\}),\|\cdot\|_{A(X,S)}$) is a quasi-normed space.
 \item[$(P2)$] Assume that $S$ has the property
  \begin{equation}\label{p2}
  \|\{a_n\}\|_S\leq C\lim_{n\to\infty} \|\{a_k\}_{k=0}^n\|_S \text{ for all } \{a_n\}\in S \text{ with } a_0\geq a_1\geq \cdots\geq 0,
  \end{equation}
  where $C>0$ is a constant depending only on $S$. Assume also that $X$ is complete and that any non-increasing sequence of non-negative numbers $\{a_n\}\in \ell_{\infty}$  belongs to $S$ if and only if $\lim_{n\to\infty}\|\{a_k\}_{k=0}^n\|_S<\infty$. Then $A(X,S,\{A_n\})$ is complete.
  \end{itemize}
\begin{definition} A scale of smoothness is a family of pairwise distinct sequence spaces $\mathbb{S}=\{S_i\}_{i\in I}$ such that $S_i\subseteq c_0$ for all $i$, all inclusions $S_i\hookrightarrow S_j$ are continuous, and the inclusion relation $\subset$ defines a total order on $\mathbb{S}$.
\end{definition}
Some examples of scales of smoothness are the Lorentz scale of smoothness $\mathbb{L}=\{\ell_{p,r}:0<r<\infty \text{ and } 0<p\leq \infty\}$ and the Lorentz-Zygmund scale of smoothness $\mathbb{L}\mathbb{Z}=\{\ell_{p,r,\gamma} : 0<r,\gamma <\infty \text{ and } 0<p\leq \infty\}$, where
 $$\ell_{p,r,\gamma}=\{\{a_n\}\in\ell_\infty: \|\{a_n\}\|_{p,r,\gamma}=\left[\sum_{n=1}^\infty n^{rp-1}(1+\log n)^{\gamma p}(a_n^*)^p\right]^\frac{1}{p}<\infty\}$$

More precisely, we have the following technical result:
\begin{lemma}\label{dos} $\mathbb{L}$ is an scale of smoothness, since the following are strict inclusions:
\begin{itemize}
\item[$(a)$] $\ell_{p,r+e}\subset \ell_{q,r}$ for all $0<p,q$ and $0<r,e<\infty$.
\item[$(b)$] $\ell_{p,r}\subset \ell_{q,r}$ for all $0<p<q$ and $0<r<\infty$.
\end{itemize}
Furthermore, $\mathbb{LZ}$ is also an scale of smoothness, since Lorentz-Zygmund sequence spaces satisfy the following strict inclusions:
\begin{itemize}
\item[$(c)$] $\ell_{p,r+e,\gamma}\subset \ell_{p,r,\alpha}$ for all $0<p\leq \infty$,$0<\gamma,\alpha<\infty$ and $0<r,e<\infty$.
\item[$(d)$] $\ell_{p,r,\gamma}\subset \ell_{q,r,\gamma}$ for all $0<p<q\leq \infty$,$0<\gamma<\infty$ and $0<r<\infty$.
\item[$(e)$] $\ell_{p,r,\gamma}\subset \ell_{p,r,\alpha}$ for all $0<p\leq \infty$, $0<r<\infty$ and $0<\alpha<\gamma <\infty$.
\end{itemize}
 \end{lemma}
\noindent \textbf{Proof. } For the case of Lorentz sequence spaces, see \cite[Lemma 1.5.2. and subsequent remarks, pp. 29-31]{stephani}. Part (a) is also proved by other means at Section \ref{generalschemes} of this paper. The proof for the Lorentz-Zygmund scale of smoothness is similar. {\hfill $\Box$}

\begin{definition} We say that the  scale of smoothness $\mathbb{S}$ is admissible with respect to the approximation scheme $(X,\{A_n\})$ if all members of $\mathbb{S}$  are admissible sequence spaces with respect to $(X,\{A_n\})$.  We say that the approximation scheme $(X,\{A_n\})$ preserves the scale of smoothness $\mathbb{S}$ if $\mathbb{S}$ is admissible and $A(X,S_i,\{A_n\})$ is strictly contained into $A(X,S_j,\{A_n\})$ whenever $S_i,S_j\in\mathbb{S}$ and $S_i\subset S_j$.
\end{definition}

\begin{definition} We say that the sequence space $S$ is rearrangement invariant if $\{a_n\}\in S \Leftrightarrow \{a_n^*\}\in S$. We say that the scale of smoothness $\mathbb{S}$ is rearrangement invariant if all members of $\mathbb{S}$  are rearrangement invariant. \end{definition}

Note that, given any sequence space $S$ satisfying $(A1)$, $(A2)$, the new sequence space given by
 \[
 S^*=\{\{a_n\}:\{a_n^*\}\in S\},
 \]
 when provided by the norm $\|\{a_n\}\|_{S^*}=\|\{a_n^*\}\|_{S}$ is rearrangement invariant. Furthermore, if $S$ is admissible with respect to the approximation scheme  $(X,\{A_n\})$, then $S^*$ is also admissible for this approximation scheme and $A(X,S,\{A_n\})=A(X,S^*,\{A_n\})$.  It is because of these observations that we assume, in all what follows, that our admissible sequence spaces and our admissible scales of smoothness are rearrangement invariant.


\section{The linear case}
\begin{theorem}Let $(X,\{A_n\})$ be a nontrivial linear approximation scheme and let us assume that $X$ is Banach. If $S_1,S_2$ are two
sequence spaces, $S_2$ is strictly contained into $S_1$, and, for all decreasing sequences $a_0\geq a_1\geq $ we have that
\begin{equation}\label{p} \{a_n\}\in S_i \text{ if and only if } \lim_{n\to\infty}\|\{a_k\}_{k=0}^n\|_{S_i}<\infty \ \ (i=1,2).
\end{equation}
  Then the norms of $A(X,S_1,\{A_n\})$ and  $A(X,S_2,\{A_n\})$ are not equivalent. In particular, if the spaces $S_i$ satisfy \eqref{p2}, then
   $A(X,S_2,\{A_n\})$ is strictly contained into $A(X,S_1,\{A_n\})$.  In other words: $(X,\{A_n\})$ preserves all scales of smoothness whose members satisfy $\eqref{p2}$ and $\eqref{p}$. Furthermore, if $X$ is Hilbert or $\dim A_n<\infty$ for all $n$ then $(X,\{A_n\})$ preserves all scales of smoothness.
\end{theorem}
\noindent \textbf{Proof.} We just need to prove that the norms of the involved spaces are not equivalent.  Now, by hypothesis, we know that there exists a sequence $\{\varepsilon_n\}\in S_1\setminus S_{2}$. Moreover, we can assume with no loss of generality that $\varepsilon_0>\varepsilon_1>\cdots>0$. Let $y_N\in A_{N+1}$ be such that $E(y_N,A_k)=\varepsilon_k$ for all $k\leq N$ (this element exists because of the linearity assumption and because $X$ is Banach, see \cite{borodin}). Obviously, we have that $\{y_N\}\subseteq A(X,S_1,\{A_n\})\cap A(X,S_2,\{A_n\})$,
\[
\lim_{N\to\infty}\|y_N\|_{A(X,S_1)}=\lim_{n\to\infty}\|\{\varepsilon_k\}_{k=0}^n\|_{S_1}<\infty
\]
and
\[
\lim_{N\to\infty}\|y_N\|_{A(X,S_2)}=\lim_{n\to\infty}\|\{\varepsilon_k\}_{k=0}^n\|_{S_2}=\infty.
\]
This proves first part of the theorem. Last claim is a direct consequence of classical Bernstein's Lethargy Theorem and its generalization to Hilbert setting by Nikolskii \cite{nikolskii}, \cite{nikolskii2} and Tjuriemskih \cite{tjuriemskih}, \cite{tjuriemskih1} (see also \cite{almiradeltoro1}){\hfill $\Box$}

\section{General approximation schemes}  \label{generalschemes}
\begin{theorem}\label{teo2} Let us assume that $(X,\{A_n\})$ is an approximation scheme such that $E(S(X)\cap A_{n+1},A_n)\geq c>0$, for infinitely many natural numbers $n\in\mathbb{N}$, and for a certain constant $c>0$. Let $S_1,S_2$ be two admissible sequence spaces verifying 
$$\lim_{N\to\infty}\frac{\|\{1,1,\cdots,1^{\text{N-th position}},0,0,\cdots\}\|_{S_1}}{\|\{1,1,\cdots,1^{\text{N-th position}},0,0,\cdots\}\|_{S_2}}=+\infty.$$ Then the norms of $A(X,S_1,\{A_n\})$ and  $A(X,S_2,\{A_n\})$ are not equivalent.  \end{theorem}
\noindent \textbf{Proof. } By hypothesis, there exists a sequence of elements $a_N\in A_N$ $(N\in N_0\subseteq \mathbb{N})$ such that $a_N\in A_{N}, \|a_N\|_X=1$, and $E(a_N,A_{N-1})>c$. Then $1\geq E(a_N,A_k)>c$ for all $k<N$ and $E(a_N,A_k)=0$ for all $k\geq N$. Let us use the notation $\mathbf{1}_N=\{1,1,\cdots,1^{\text{N-th position}},0,0,\cdots\}$. Then
$$\|\{E(a_N,A_k)\}_{k=0}^\infty\|_{S_1}\geq c\|\mathbf{1}_N\|_{S_1}\text{ and  } \|\{E(a_N,A_k)\}_{k=0}^\infty\|_{S_2}\leq \|\mathbf{1}_N\|_{S_2},$$
so that
\[
\frac{\|a_N\|_{A(X,S_1)}}{\|a_N\|_{A(X,S_2)}}=\frac{\|\{E(a_N,A_k)\}_{k=0}^\infty\|_{S_1}}{\|\{E(a_N,A_k)\}_{k=0}^\infty\|_{S_2}}\geq c\frac{\|\mathbf{1}_N\|_{S_1}}{\|\mathbf{1}_N\|_{S_2}}\to\infty \text{ (for } N\to\infty\text{)}.
\]
This proves the theorem.  {\hfill $\Box$}

\noindent \textbf{Remark.} Theorem \ref{teo2} can be used for the Lorentz sequence spaces $\ell_{p,r}$. To prove this it is enough to take into account that, for all $\alpha>-1$,
\[
\lim_{N\to\infty}\frac{\sum_{k=1}^Nk^{\alpha}}{N^{\alpha+1}}=\frac{1}{\alpha+1},
\]
(see \cite[Part I, Problem 71]{polyaszego}) so that
\begin{eqnarray*}
\lim_{N\to\infty}\frac{\|\mathbf{1}_N\|_{\ell_{q,r_2}}}{\|\mathbf{1}_N\|_{\ell_{p,r_1}}} &=& \lim_{N\to\infty}\frac{(\sum_{k=1}^N k^{r_2q-1})^{\frac{1}{q}}}{(\sum_{k=1}^Nk^{r_1p-1})^{\frac{1}{p}}} \\
&=&  \lim_{N\to\infty}\frac{(N^{r_2q}/r_2q)^{\frac{1}{q}}}{(N^{r_1p}/r_1p)^{\frac{1}{p}}} \\
&=&  \frac{(r_1p)^{\frac{1}{p}}}{(r_2q)^{\frac{1}{q}}} \lim_{N\to\infty}N^{r_2-r_1} \\
&=& \left\{\begin{array}{lllll}
+\infty &  &  \text{ for }r_1<r_2 \\
0 &  &  \text{ for }r_1>r_2 \\
\frac{(r_1p)^{\frac{1}{p}}}{(r_2q)^{\frac{1}{q}}} &  &  \text{ for }r_1=r_2 \\
\end{array}
\right.
\end{eqnarray*}
In particular, these computations show that inclusion $\ell_{p,r+e}\subset \ell_{q,r}$ is strict for all $r,e>0$ and all $p,q>0$. Unfortunately, these computations are not useful for the important case $r_1=r_2$.
In order to deal with this case, so that we can prove that our approximation scheme preserves the Lorentz smoothness scale, we need to restrict a little bit more the class of approximation schemes we are considering. Concretely, we can use the following result by Brudnyi:
\begin{theorem}[Brudnyi, see Theorem 4.5.12 and Remark 4.5.13 in \cite{brundykrugljak}]\label{teobrud} Let us assume that $(X,\{A_n\})$ is an approximation scheme such that $X$ is Banach and  $E(S(X)\cap A_{n+1},A_n)\geq c>0$ for all $n\in\mathbb{N}$, and for a certain constant $c>0$. Let us also assume that $A_n+A_m\subseteq A-{n+m}$ for all $n,m\in\mathbb{N}$ and let $\{\varepsilon_n\}$ be any convex non-increasing sequence of positive real numbers belonging to $c_0$. Then:
\begin{itemize}
\item[$(i)$] There exists an element $x\in X$ such that $E(x,A_n)\geq \varepsilon_n$ for all $n\in\mathbb{N}$ and $\lim\inf \frac{E(x,A_n)}{\varepsilon_n}<\infty$.
\item[$(ii)$] If, furthermore, $\sup_{n\in\mathbb{N}}{\frac{\varepsilon_n}{\varepsilon_{2n}}}<\infty$ then there exists $x\in X$ such that $E(x,A_n)\approx \varepsilon_n$
\end{itemize}
  \end{theorem}
\begin{corollary}\label{corbrud} If the approximation scheme $(X,\{A_n\})$ verifies the hypotheses of Theorem \ref{teobrud}, then it preserves the Lorentz scale of smoothness.
\end{corollary}
\noindent \textbf{Proof. } We only need to deal with the case $r_1=r_2=r$ and $p<q$. Take $\varepsilon_n=n^{-r}(1+\log_2n)^{-\frac{1}{q}}$. This sequence belongs to $\ell_{r,q}\setminus\ell_{r,p}$ and it is convex. Furthermore,
\[
\sup_{n\in\mathbb{N}}{\frac{n^{-r}(1+\log_2n)^{-\frac{1}{q}}}{(2n)^{-r}(1+\log_2(2n))^{-\frac{1}{q}}}}= \sup_{n\in\mathbb{N}}{\frac{1}{2^{-r}}\frac{n^{-r}(1+\log_2n)^{-\frac{1}{q}}}{n^{-r}(2+\log_2n)^{-\frac{1}{q}}}}<\infty
\]
Thus, if we apply part $(ii)$ of Brudnyi's Theorem, we have that there exists $x\in X$ such that  $E(x,A_n)\approx \varepsilon_n$. Hence $x\in A_q^r(X,\{A_n\})\setminus A_p^r(X,\{A_n\})$. This ends the proof. {\hfill $\Box$}

An interesting example of approximation scheme verifying Brudnyi's conditions is the following one: Take $X$ a quasi-Banach space such that there exists an infinite sequence of linear projections $P_n:X\to X$ satisfying $\rank(P_n)=n$ for all $n\geq 1$ and $\sup_{n\geq 1}\|P_n\|=C<\infty$ (e.g. this condition holds for $X$ whenever $X$ contains a complemented subspace $Y$ such that $Y$ has a Schauder basis). Then $(X,\{\Sigma_n\})$ satisfies  the hypotheses of Theorem \ref{teobrud}, where $\Sigma_n=\{T\in B(X,X):\rank(T)<n\}$. Of course, in this case, $E(T,\Sigma_n)=a_n(T)$ is the $n$-th approximation number of the operator $T$. To prove this claim, it is enough to take into account that, if $H_n=P_n(X)$ and we define $Q_n:X\to H_n$ by $Q_n(x)=P_n(x)$ and $i_n:H_n\to X$ denotes the inclusion map, then $1_{H_n}=Q_nP_ni_n$, $\|Q_n\|=\|P_n\|\leq C$, $\|i_n\|=1$, so that:
\[
1=a_n(1_{H_n})\leq \|Q_n\|a_n(P_n)\|i_n\|\leq Ca_n(P_n) \ \text{for all } n=1,2,\cdots.
\]
Hence $a_n(\frac{P_n}{\|P_n\|})=\frac{1}{\|P_n\|}a_n(P_n)\geq \frac{1}{C}a_n(P_n)\geq c=1/C^2>0$, $n=1,2,\cdots$, so that $\frac{P_n}{\|P_n\|}\in S(B(X,X))\cap \Sigma_{n+1}$ and $E(\frac{P_n}{\|P_n\|},\Sigma_n)=a_n(\frac{P_n}{\|P_n\|})>c$.

\begin{corollary}\label{corcor} Let us assume that $X$ is Banach and the approximation scheme $(X,\{A_n\})$ satisfies Shapiro's theorem. Then there exists an increasing  sequence of natural numbers $\{m(n)\}$ such that $(X,\{A_{m(n)}\})$ preserves the Lorentz scale of smoothness. \end{corollary}

\noindent \textbf{Proof. } By hypothesis, $(X,\{A_n\})$ satisfies Shapiro's theorem so that $E(S(X),A_n)=1$ for all $n$ (see \cite[Corollary 3.6]{almira}). It follows that we can use the density of $\bigcup_{n\in\mathbb{N}}A_n$ in $X$ to prove that, for a certain increasing sequence of natural numbers $m(n)$,   $B_n=A_{m(n)}$ satisfies $B_n+B_m\subseteq B_{n+m}$ and $E(S(X)\cap B_{n+1},B_n)\geq \frac{1}{2}$ for all $n$, so that $(X,\{B_n\})$ is an approximation scheme verifying the hypotheses of Theorem \ref{teobrud}. {\hfill $\Box$}.

\noindent \textbf{Remark.} There are, of course, non-linear approximation schemes verifying stronger results than the general results given by Theorem \ref{teo2} and Corollaries \ref{corbrud}, \ref{corcor}. For example, Starovoitov \cite{starovoitov} has proved that approximation of continuous functions by rational functions $R_n:=\{p(t)/q(t): p,q\in \Pi_n \text{ and } q(t)\neq 0 \text{ for all } t\in[a,b]\}$, satisfies Bernstein's Lethargy Theorem, so that the non-linear approximation scheme $(C[a,b],\{R_n\})$ preserves all scales of smoothness.


\section{The case of approximation by finite rank operators}
 In this section we will study, for arbitrary infinite dimensional Banach spaces $X,Y$, the inclusion relations between the ideal of operators  $\mathcal{L}_{p,r}^{(a)}(X,Y)=\{T\in B(X,Y): (a_n(T))\in \ell_{p,r}\}$, where $a_n(T)=\inf_{S\in B(X,Y), \rank(S)<n}\|T-S\|$ denotes the $n$-th approximation number of the operator $T$. These ideals have been the subject of many papers and monographs. In particular, if we use the standard notation $\mathcal{L}_{p,r}^{(a)}=\bigcup_{X,Y} \mathcal{L}_{p,r}^{(a)}(X,Y)$, and we take into account the inclusions of Lorentz sequence spaces (see Lemma \ref{dos} below), then it is clear that  $\mathcal{L}_{p,r+e}^{(a)}\subseteq \mathcal{L}_{q,r}^{(a)}$ for all $0<p,q$ and $0<r,e<\infty$ and  $\mathcal{L}_{p,r}^{(a)}\subseteq \mathcal{L}_{q,r}^{(a)}$ for all $0<p<q$ and $0<r<\infty$. What is more, it is well known that, via the computation of approximation numbers of the diagonal operators $D:\ell_p\to\ell_p$, it is possible to show that these inclusions are strict. But, to the knowledge of the authors of this paper, strict inclusions of the type $\mathcal{L}_{p,r+e}^{(a)}(X,Y)\subset \mathcal{L}_{p,r}^{(a)}(X,Y)$ or $\mathcal{L}_{p,r}^{(a)}(X,Y)\subset \mathcal{L}_{q,r}^{(a)}(X,Y)$ have not been yet proved  for arbitrary Banach spaces $X,Y$. This is what we make, with the help of a very strong result by Oikhberg \cite{oikhberg}, in this section. Concretely, we prove that for all $X,Y$ infinite dimensional Banach spaces, approximation by finite rank bounded linear operators  $T:X\to Y$ preserves the Lorentz scale of smoothness. A similar result also holds true for Lorentz-Zygmund scale of smoothness.  Some related results appear in
\cite{tita_equiv_quasinorm}, where equivalence of several distinct norms for some operator ideals (defined as approximation spaces with the help of symmetric norming functions) are proved.

\begin{lemma}\label{uno} Assume that $\{a_n\}\in \ell_{p,r}$ and $0<C<\infty$. Then $(a_{[n/C]})\in \ell_{p,r}$. \end{lemma}
\noindent \textbf{Proof. } By definition, $\{a_n\}\in \ell_{p,r}$ if and only if $\{a_n^*\}\in \ell_{p,r}$, so that we can assume with no loss of generality, that $\{a_n\}$ is non-increasing. Then   $(a_{[n/C]})$ is also non-increasing and, if we take into account that  $[n/C]$ can repeat its value at most  $[C]+1$ times, we get:
\begin{itemize}
\item Case $p\geq\frac{1}{r}$:
\begin{eqnarray*}
\sum_{n=1}^\infty n^{rp-1}(a_{[n/C]})^p &\leq& \sum_{n=1}^\infty \left(C([n/C]+1)\right)^{rp-1}(a_{[n/C]})^p \\
&\leq&  \sum_{m=1}^\infty ([C]+1)\left(C(m+1)\right)^{rp-1}a_{m}^p \\
&\leq&  ([C]+1)C^{rp-1}2^{p r-1} \sum_{m=1}^\infty m^{rp-1}a_{m}^p <\infty,
\end{eqnarray*}
since $\sup_{m\geq 1} \frac{(m+1)^{rp-1}}{m^{rp-1}}=2^{rp-1}$.
\item Case $p<\frac{1}{r}$: \begin{eqnarray*}
\sum_{n=1}^\infty n^{rp-1}(a_{[n/C]})^p &\leq& \sum_{n=1}^\infty \left(C[n/C]\right)^{rp-1}(a_{[n/C]})^p \\
&\leq&  ([C]+1)C^{rp-1}\sum_{m=1}^\infty m^{rp-1}a_{m}^p. \end{eqnarray*} \end{itemize}
This ends the proof.  {\hfill $\Box$}

Indeed, a more general result can be proved.
\begin{lemma}\label{tres} Let $\alpha,C>1$ and let $S$ be an admissible sequence space with respect to approximation schemes satisfying  $K(n)=\alpha n$. Let $\{a_n\}$ be a non-increasing sequence of non-negative real numbers. Then  $\{a_n\}\in S \Leftrightarrow (a_{[n/C]})\in S$. In particular, this equivalence is verified by all Lorentz and Lorentz-Zygmund sequence spaces, since they are admissible for $K(n)=2n$.
\end{lemma}
\noindent \textbf{Proof. } Let us assume that $(a_{[n/C]})\in S$. The inequality $a_{[n/C]}\geq a_n$ for all $n$, in conjunction with the fact that $S$ is a solid, implies that $\{a_n\}\in S$. Let us prove the other implication. We use that $\alpha>1$ to claim that, for a certain natural number $k\geq 1$, $\alpha^k>[C]+1$. Let us assume that $\{a_n\}$ is a non-increasing sequence of non-negative real numbers and let $\{b_n\}$ be the sequence given by $b_n=a_{[n/C]}$, $n=1,2,...$. We want to show that $\{b_n\}\in S$.  Now, $b_{\alpha^kn}=a_{\left[\frac{\alpha^kn}{C}\right]}\leq a_n$ for all $n$ (since $\{a_n\}$ is non-increasing and $\alpha^kn/C\geq ([C]+1)n/C\geq n$). Hence $\{b_{\alpha^kn}\}\in S$. Now, $\{b_{\alpha^{k-1}n}\}$ is non-increasing, so that admissibility of $S$ for $K(n)=\alpha n$ implies that $\{b_{\alpha^{k-1}n}\} \in S$. A repetition of this argument $k$ times gives us $\{b_n\}\in S$. This ends the proof. {\hfill $\Box$}

\begin{theorem}\label{teo} For all $X,Y$ infinite dimensional Banach spaces, approximation by finite rank bounded linear operators  $T:X\to Y$ preserves all its admissible scales of smoothness. In particular, it preserves the Lorentz and the Lorentz-Zygmund scales of smoothness.
\end{theorem}
\noindent \textbf{Proof. } Let $\mathbb{S}$ be an admissible scale of smoothness for approximation by finite rank bounded linear operators and let us assume that $S_1,S_2\in\mathbb{S}$, $S_1\subset S_2$. Let $(\varepsilon_n)\in S_2\setminus S_1$.  We can assume with no loss of generality that $(\varepsilon_n)$ is non-increasing and converges to zero, since our sequence spaces are assumed to be rearrangement invariant. Now, Theorem 1.1 from \cite{oikhberg} guarantees that there exists an operator $T\in B(X,Y)$ such that $3\varepsilon_{[n/6]}\geq a_n(T)\geq \varepsilon_{n}/9$ for all $n=1,2,\cdots$.   Hence  Lemma \ref{tres} implies that $(a_n(T))\in S_2\setminus S_1$. This ends the proof. {\hfill $\Box$}


\section{Real interpolation spaces}
As it is well known, central theorems in approximation theory state a strong relation between approximation spaces and the interpolation spaces obtained by the use of the $K$-functional of Petree
$K(x,t,X,Y)=\inf_{y\in Y}\|x-y\|_X+t\|y\|_Y$,
\[
(X,Y)_{\theta,q}=\{x\in X: \|x\|_{\theta,q}=\|t^{-(\theta+\frac{1}{q})}K(x,t,X,Y)\|_{L^q(0,\infty)}<\infty\}, \ \ (0<\theta<1, 0<q\leq \infty).
\]
In particular, if we use that $Y$ is continuously embedded into $X$ and $K(x,t,X,Y)$ is a monotone function of $t$, it is not difficult to prove that the norm of $(X,Y)_{\theta,q}$ is equivalent to the norm $\rho_{\theta,q}(x)=\|\{2^{k\theta}K(x,\frac{1}{2^k},X,Y)\}_{k=0}^{\infty}\|_{\ell_q}$ and, with this new norm, interpolation spaces and approximation spaces are quite similar objects: we just replace the errors of best approximation by the evaluation of Petree's $K$-functional at the points of a certain decreasing sequence of positive real numbers. It is then natural to ask for a version of Bernstein's lethargy theorem in terms of $K$-funcionals. This result already exists and was obtained by Krugljak \cite[Theorems 4.5.7 and 4.5.10]{brundykrugljak}. Concretely, he proved that, given a couple of quasi-Banach spaces $(X_0,X_1)$, the following are equivalent claims:
\begin{itemize}
\item[(K1)] For each continuous concave function $\varphi:(0,1]\to [0,\infty]$ such that $\lim_{t\to 0}\varphi(t)=0$ there exists an element $x\in X_0+X_1$ such that $K(x,t,X_0,X_1)\approx \varphi$.
\item[(K2)] There exists $x\in X_0+X_1$ such that $\int_0^t K(x,s,X_0,X_1)ds\leq \gamma K(x,t,X_0,X_1)$ for all $t\in (0,1]$ and a certain $\gamma>0$.
\end{itemize}
An easy consequence  of Krugljak's theorem is the following
\begin{corollary} Let us assume that $(X,Y)$ is a couple, and $Y\subset X$. Let us assume that condition $(K2)$ is satisfied for the couple $(X,Y)$ and let $(\theta,q),(\theta^*,q^*)$ be two distinct points of $(0,1)\times (0,\infty]$. Then $(X,Y)_{\theta,q}\neq (X,Y)_{\theta^*,q^*}$. Furthermore, all these spaces are strictly contained into $X$.
\end{corollary}

In this section we state and prove a version of Shapiro's theorem in terms of $K$-functionals and we use it to state a general condition for the strict inclusion of $(X,Y)_{\theta,q}$ into $X$.

\begin{theorem}\label{kfunc}Let $Y$ be a quasi-semi-normed subspace of the quasi-Banach space $X$ and let us consider the $K$-funcional $K(x,t)=\inf_{y\in Y}[\|x-y\|_X+t\|y\|_Y]$. The following are equivalent claims:
\begin{itemize}
\item[$(a)$] There exists  $c>0$ such that, for all $t\in (0,1]$,  $\sup_{\|x\|_X=1}K(x,t)>c$.
\item [$(b)$] There exists  $c>0$ such that, for every strictly decreasing sequence of positive numbers $\{t_n\}\in c_0$, we have that  $\sup_{\|x\|_X=1}K(x,t_n)>c$, $n=0,1,2,\cdots$.
\item [$(c)$] There exists  $c>0$ and a strictly decreasing sequence of positive numbers $\{t_n\}\in c_0$ such that  $\sup_{\|x\|_X=1}K(x,t_n)>c$, $n=0,1,2,\cdots$.
\item[$(d)$] For all $(b_n)\subset [0,\infty)$ such that $\lim_{n\to\infty}b_n=\infty$ and all  strictly decreasing sequence of positive numbers $\{t_n\}\in c_0$, we have that $A(\{b_n\},\{t_n\})=\{x\in X: \sup_{n\geq 1}b_nK(x,t_n)<\infty\}$ is a proper subset of $X$.
\end{itemize}
In particular, when $(a)$ holds true, the interpolation spaces $(X,Y)_{\theta,q}$ are proper subspaces of $X$.
\end{theorem}
\noindent \textbf{Proof. } The monotonicity of $K(x,t)$ on $(0,1]$ gives the equivalences $(a) \Leftrightarrow (b) \Leftrightarrow (c)$.
To prove the other equivalences we need firstly note that $A(\{b_n\},\{t_n\})$ is a quasi-Banach space continuously embedded into $X$.
Let us now show that $(b) \Leftrightarrow (d)$.

\noindent $(b)\Rightarrow (d)$ Let, for each $n\in \mathbb{N}$, $x_n\in S(X)$ be such that $K(x_n,t_n)>c$. Then
\[
\|x_n\|_{A(\{b_n\},\{t_n\})}=\sup_{m\geq 1}b_mK(x_n,t_m)\geq b_n K(x_n,t_n)>b_nc\to\infty.
\]
This shows that the norms of $A(\{b_n\},\{t_n\})$ and $X$ are not equivalent, so that $A(\{b_n\},\{t_n\})\neq X$.  This ends the proof.

\noindent $(d)\Rightarrow (b)$ We show that the negation of $(b)$ implies the negation of $(d)$.  Let $\{t_n\}\in c_0$ be a decreasing sequence and let us consider the values  $d_n= \sup_{\|x\|_X=1}K(x,t_n)$. Then $(d_n)$ is a non-increasing sequence, so that $(b)$ fails for $\{t_n\}$ if and only if $d_n\to 0$. Assume this is the case. Let $b_n=\frac{1}{d_n}$ and let $x\in X$. Then
\[
\sup_{n\geq 1}b_nK(x,t_n)=\sup_{n\geq 1}b_n\|x\|_XK(\frac{x}{\|x\|_X},t_n)\leq \sup_{n\geq 1}b_n\|x\|_Xd_n\leq \|x\|_X<\infty.
\]
Hence $A(\{b_n\},\{t_n\})=X$ and $(d)$ also fails for $\{t_n\}$. This ends the proof of the equivalences.

Let us now assume that condition $(a)$ holds true. In order to prove that $(X,Y)_{\theta,q}$ is a proper subspace of $X$, we use that $(X,Y)_{\theta,q}\subseteq (X,Y)_{\theta,\infty} = A(\{2^{n\theta}\},\{2^{-n}\}) \subsetneqq X$ {\hfill $\Box$}

Sometimes condition $(a)$ of Theorem \ref{kfunc}  can be checked directly. For example, $(a)$ holds true as soon as $\overline{Y}^X$ is properly contained into $X$, since in that case there exists $x\in X$ with $0<c=E(x,Y)\leq K(x,t)$ (all $t>0$). If we assume that $Y$ is a dense subspace of $X$, things can be more complicated, but there are also cases where a simple computation makes the work. For example, for $X=C[0,1]$ and $Y=C^{(1)}[0,1]$ it is easy to find, for each $n\geq 2$, and for $0<a<b<1$, $f(t)\in C[0,1]$ such that $f(a)=-1$, $f(b)=1$. Hence, if $g\in C^{(1)}[0,1]$ and $\|f-g\|_{\infty}<1/2$ then $g(a)<-1/2$, $g(b)>1/2$, and the Mean Value Theorem guarantees that $\|g\|_{C^{(1)}}\geq \frac{g(b)-g(a)}{b-a}\geq 1/(b-a)$. Thus, for $t< 1/2$ we have that $K(f,t)\geq t\frac{1}{b-a}$. Taking $b-a\leq t$  we get $K(f,t)\geq 1/2$. Of course, these computations are just a particular case of the following general situation:

\begin{proposition}Let $X$, $Y$ be such that $Y$ is a dense subspace of $X$ and there exists an strictly decreasing function $\phi:(0,1]\to\mathbb{R}^+$ and a constant $c\in (0,1)$ such that $\lim_{t\to 0}\phi(t)=+\infty$ and, for each $\varepsilon>0$ there exists $x_{\varepsilon}\in S(X)$ such that $y\in Y$ and $\|x_{\varepsilon}-y\|_X<c$ implies $\|y\|_Y\geq \phi(\varepsilon)$. Then the $K$-functional associated to the pair $(X,Y)$ satisfies condition $(c)$ of Theorem \ref{kfunc}
\end{proposition}
\noindent \textbf{Proof. } Take $t_n=\frac{1}{\phi(1/n)}$ and $z_n=x_{\frac{1}{n}}$. Then $K(z_n,t_n)\geq c$ since otherwise we would have $\|z_n-y\|_X+t_n\|y\|_Y<c$ for a certain $y\in Y$, which is impossible since $\|z_n-y\|_X\leq c$ implies $t_n\|y\|_Y\geq t_n\phi(1/n)=1>c$. This ends the proof. {\hfill $\Box$}

Of course, in many concrete cases, $K(x,t)$ is known (or, at least, an equivalent function $w(x,t)\approx K(x,t)$ is known). In these cases part $(a)$ can also be checked by direct computations and, usually, this checking can be made with some easy estimations.

On the other hand, it would be nice to know if condition $(a)$ of Theorem \ref{kfunc} follows just from the strictness of the inclusions $(X,Y)_{\theta,q}\subset X$, since a similar result is known for approximation spaces. Concretely, in \cite{almira} it is proved that approximation schemes $(X,\{A_n\})$ verifying $K(n)=cn$ satisfy Shapiro's theorem  if and only if the approximation space $A_q^r(X,\{A_n\})$ is properly contained into $X$ for some choice of parameters $q,r$.  Finally,  another main open question (still unsolved) is to know if condition $(K2)$ above can be relaxed in order to guarantee that the natural inclusions between interpolation spaces are strict.

\section{Acknowledgements} I would like to express my gratitude to Prof. N. Ti\c{t}a from Brasov  for many fruitful discussions on the subject of this paper.


\bigskip

\footnotesize{J. M. Almira

Departamento de Matem\'{a}ticas. Universidad de Ja\'{e}n.

E.P.S. Linares,  C/Alfonso X el Sabio, 28

23700 Linares (Ja\'{e}n) Spain

Email: jmalmira@ujaen.es}

\end{document}